\documentclass[a4paper, fontsize=12pt]{article}

\renewenvironment{abstract}
               {\list{}{\rightmargin\leftmargin}%
                \item[\hspace*{1cm}\small\textbf{Abstract ---}]\relax}
               {\endlist}

\usepackage[top=1in, bottom=1in, left=1in, right=1in]{geometry}
\usepackage[fleqn]{amsmath}
\usepackage{amssymb}
\usepackage{enumerate}
\usepackage{amsthm}
\usepackage{appendix}
\usepackage[russian,english]{babel}
\usepackage{subfigure}
\usepackage{float}
\usepackage[pdftex]{graphicx}
\newtheorem{Theorem}{Theorem}

\newtheorem{Agreement}[Theorem]{Agreement}
\newtheorem{Definition}[Theorem]{Definition}
\newtheorem{Remark}[Theorem]{Remark}

\newtheorem{Example}[Theorem]{Example}

\begin{document}

\title{\bf Hyperreal delta functions as a new general tool for modeling physical states with infinitely high densities}

\author{
        Marcoen J.T.F. Cabbolet\footnote{e-mail: Marcoen.Cabbolet@vub.ac.be}\\
        \small{\textit{Center for Logic and Philosophy of Science, Vrije Universiteit Brussel}}\\
        }
\date{}

\maketitle

\begin{abstract}
This paper introduces the expanded real numbers as an ordered subring of the hyperreal number field that does not contain any infinitesimals, and defines the set of all integrable functions from the real numbers to the expanded real numbers. This allows to identify the Dirac delta with a special hyperreal-valued function of a real variable: the Dirac delta function thus defined is a general tool, applicable for the mathematical modeling of physical systems in which infinitely high densities occur.
\end{abstract}
\section{Introduction}
In physics we are interested in modeling a state of a system in which a physical quantity is distributed over space. Let's limit the discussion to one-dimensional physical systems for the sake of simplicity, and suppose then that a real amount $\xi$ of a physical quantity (e.g. mass) is distributed over the space $\mathbb{R}$: this state can then, in general, be modeled by a function $f: \mathbb{R}\rightarrow \mathbb{R}$, satisfying $\int^\infty_{-\infty} f(x)dx = \xi$, with the function value $f(x) \in \mathbb{R}$ at a point $x \in \mathbb{R}$ representing the \emph{density}. However, troubles arise when the physical quantity is distributed over isolated points $x_1, x_2, \ldots, x_n$ in the space $\mathbb{R}$: we have this situation, for example, when we consider the distribution of mass in a system made up of $n$ point particles. In such cases \emph{infinitely high densities} occur at those isolated points $x_j \in \mathbb{R}$: no real function $f$ exists that can model such a state.

To model states in which infinitely high densities occur, functions $f$ are required that can have an infinitely large value at isolated points in space, but such that $f$ can also be added to a (piecewise) smooth function $g$ representing a distribution of an amount $\chi$ that physical quantity over a region of space---a sum $f+g$ then represents a distribution of an amount $\xi+\chi$ of that physical quantity over space. Analyzing, the problem is thus that functions $f$ are required that must satisfy the following two conditions:
\begin{equation}\label{eq:0}
\begin{array}{ll}
{\rm (i)} & f: \mathbb{R} \rightarrow R\\
{\rm (ii)} & \left\{
\begin{array}{l}
f(x) = 0 \Leftrightarrow x \neq r \\
\int_{-\infty}^{+\infty}f(x)dx =\xi
\end{array}
\right.
\end{array}
\end{equation}
where $r, \xi \in \mathbb{R}$ and $R$ is a number ring containing the reals: $R \supset \mathbb{R}$.

The purpose of this paper is to solve this problem: below we introduce special hyperreal functions of a real variable---`special' because they satisfy both clauses of Eq. (\ref{eq:0})---that can be applied in general for the mathematical modeling of physical states in which infinitely high densities occur.

\section{Motivation for introducing new definitions}
Our attention is drawn to the Dirac delta. At its original introduction it was not exactly defined, but rather heuristically characterized as an object---denoted by the symbol $\delta$---for which
\begin{equation}\label{eq:1}
\left\{
\begin{array}{l}
\delta(x) = 0 \Leftrightarrow x \neq 0 \\
\int_{-\infty}^{+\infty}\delta(x)dx =1
\end{array}
\right.
\end{equation}
where $x$ is a real variable \cite{Dirac}. This corresponds to clause (ii) of Eq. (\ref{eq:0}) for $r = 0, \xi = 1$. There is no function $f: \mathbb{R} \rightarrow \mathbb{R}$ that has these properties. In fact, Von Neumann dismissed the idea of the Dirac delta as ``fiction'' \cite{VonNeumann}.

Since then, however, various objects that capture the idea of the Dirac delta have been rigorously defined within the framework of standard analysis. An example is the \emph{Dirac delta distribution}, a linear functional $\delta: \mathcal{D}(\mathbb{R}) \rightarrow \mathbb{R}$ on the space of test functions on $\mathbb{R}$ defined by $\langle \delta, f \rangle := f(0)$ for all $f \in \mathcal{D}(\mathbb{R})$ \cite{Schwartz}. Its definition is often written as
\begin{equation}\label{eq:abusive}
\langle \delta, f \rangle = \int_\mathbb{R} f(x)\delta(x)dx = f(0)
\end{equation}
but the objection is then that this is an abuse of notation since the term `$\delta(x)$' in the integrand does not refer to an existing object. To meet that objection, the term $\delta(x)$ can be viewed as a weak limit of a sequence of real functions $g_n:\mathbb{R} \rightarrow \mathbb{R}$, e.g. the functions $g_n$ given, for positive integers $n$, by $g_n(x) = n$ for $x \in (-\frac{1}{2n}, \frac{1}{2n})$ and $g_n(x) = 0$ else; we then get
\begin{equation}\label{eq:2}
\int_\mathbb{R} f(x)\delta(x)dx  := \int_\mathbb{R} \lim_{n\rightarrow \infty}f(x)g_n(x)dx
\end{equation}
Still, this weak limit $\delta$ does not exist in the space of all functions from $\mathbb{R}$ to $\mathbb{R}$.

Another standard object that formalizes the idea of the Dirac delta is the \emph{Dirac measure}, a real function $\delta_0: \mathcal{B}(\mathbb{R}) \rightarrow \mathbb{R}$ on the $\sigma$-algebra of Borel sets of $\mathbb{R}$ defined, for $X \in \mathcal{B}(\mathbb{R})$, by $\delta_0(X) = 1$ if $0 \in X$ and $\delta_0(X) = 0$ else \cite{Rudin}. The Lebesgue integral of a (measurable) function $f: \mathbb{R} \rightarrow \mathbb{R}$ against the measure $\delta_0$ then gives
\begin{equation}
\int_\mathbb{R} f(x)\delta_0(dx)  = f(0)
\end{equation}
So, the Dirac measure corresponds to a Dirac delta distribution as defined above in a natural way, but it has to be emphasized that the Dirac measure has no Radon derivative---that is, the Dirac measure does not correspond to an ordinary function $\delta$ on the reals such that $\int_\mathbb{R} f(x)\delta_0(dx) = \int_\mathbb{R} f(x)\delta(x)dx$.

Other objects equivalent to the above distribution and measure have been introduced, but the point is this: although they all capture the idea of the Dirac delta, none of them actually defines a function satisfying Eq. (\ref{eq:1}). That is, the objects that have been defined in the framework of standard analysis to formalize the idea of the Dirac delta are \underline{useful} for doing calculations, but are \underline{useless} for our present purpose---which is to model the state of a system in which a physical quantity is distributed over isolated points in space by a function on $\mathbb{R}^n$.

Another development is that the real number field has been extended to the hyperreal number field $^*\mathbb{R}$ \cite{Robinson}. Todorov has shown that there exists a nonstandard (hyperreal) function $^*\delta: \ ^*\mathbb{R} \rightarrow \ ^*\mathbb{R}$ such that
\begin{equation}
\int_{^*\mathbb{R}} f(x)^*\delta(x)dx = f(0)
\end{equation}
for any continuous function $f$ on $\mathbb{R}$ \cite{Todorov}. A concrete example is the function $^*\delta: \ ^*\mathbb{R} \rightarrow \ ^*\mathbb{R}$ for which
\begin{equation}\label{eq:1a}
\left\{
\begin{array}{l}
^*\delta(x) = 0 \Leftrightarrow x \not \in [-\frac{dx}{2}, \frac{dx}{2}] \\
^*\delta(x) = \frac{1}{dx}\Leftrightarrow x \in [-\frac{dx}{2}, \frac{dx}{2}]
\end{array}
\right.
\end{equation}
where `$dx$' is the infinitesimal hyperreal number $\omega^{-1}$ \cite{Dannon}. So, this is an object defined in the framework of nonstandard analysis that captures the idea of the Dirac delta of Eq. (\ref{eq:1}), but again this object is not useful for our present purposes: the point is that the function $^*\delta$ violates condition (i) of Eq. (\ref{eq:0}) because it is a function of a hyperreal variable, not of a real variable. This is objectionable, because in physics we want to model a spatial dimension with the real numbers, not with the hyperreal numbers.

Thus speaking, we cannot but conclude that the existing objects that formalize the idea of the Dirac delta are \emph{not suitable} for our present purposes: we will, thus, have to develop a new object.

\section{Expanded real delta functions of a real variable}
For our present purpose we need hyperreal numbers, but we do not need \emph{all} of the hyperreal number field. So first of all we apply Ockham's razor and we define the ordered ring of the expanded reals as the part of the hyperreals that we need:

\begin{Definition}\label{def:1}\ The \textbf{ordered ring of expanded real numbers} is the subring $(^*_+\mathbb{R}, + \ , \ \cdot \ , >)$ of the hyperreal number field $(^*\mathbb{R}, + \ , \ \cdot \ , >)$ given by
\begin{equation}\label{eq:4}
^*_+\mathbb{R} = \{\xi \in \ ^*\mathbb{R} \ | \ \xi = a_1\omega^{p_1}+\ldots+a_n\omega^{p_n}, n \in \mathbb{N}^+,  p_1 > \ldots > p_n \geq 0, a_j \in \mathbb{R}\}
\end{equation}
That is, the set $^*_+\mathbb{R}$ contains the real numbers and, as indicated by the subscript `+' on the left, those hyperreal numbers $\xi \not \in \mathbb{R}$ in which a finite number of \underline{positive} powers of the infinitely big hyperreal number $\omega$ with $|\omega| = \infty$  occur. As a set, expanded real numbers can, for \underline{nonzero} $a_j$, be defined according to
\begin{equation}\label{eq:7}
a_1\omega^0  :=  a_1
\end{equation}
\begin{equation}\label{eq:8}
a_2\omega^n := \{ \langle a_2\omega^{n-1},0\rangle, \langle a_2\omega^{n-1},1\rangle, \langle a_2\omega^{n-1},2\rangle, \ldots \}
\end{equation}
\begin{equation}\label{eq:8a}
\sum_{j=1}^n a_j\omega^{p_j} := a_1\omega^{p_1} \cup a_2\omega^{p_2} \cup \ldots \cup a_n\omega^{p_n}
\end{equation}
with $a_1$ in Eq. (\ref{eq:7}) to be represented by a Dedekind cut, $n \geq 1$ in Eq. (\ref{eq:8}), and $p_1 > p_2 > \ldots > p_n \geq 0$ in Eq. (\ref{eq:8a}). Note that then $^*_+\mathbb{R} \supset\mathbb{R}$, with an inequality $a_1 \omega + a_2 \neq a_1 \omega$ being an inequality of sets.
\hfill $\Box$
\end{Definition}

\begin{Agreement}\rm\ We will henceforth take the notation $\sum_{j=1}^n a_j\omega^{p_j}$ for an expanded real number $x \in\ ^*_+\mathbb{R}$ to \emph{imply} that $p_1 > p_2 > \ldots > p_n \geq 0$ as in Def. \ref{def:1} above.
\hfill $\Box$
\end{Agreement}

\begin{Remark}\rm\ It has to be taken that Eqs. (\ref{eq:7})-(\ref{eq:8a}) give a \emph{set-theoretical representation} of a \textbf{single} expanded real number. In particular, Eq. (\ref{eq:8a}) means that `$\sum_{j=1}^n a_j\omega^{p_j}$' is a \textbf{notation} for the set $a_1\omega^{p_1} \cup a_2\omega^{p_2} \cup \ldots \cup a_n\omega^{p_n}$. This does \textbf{absolutely not} mean that the binary operation `addition` on the ring of expanded real numbers is identified with taking the union of two sets, as in
\begin{equation}
\alpha + \beta \equiv \alpha \cup \beta
\end{equation}
for any $\alpha,\beta\in\ _+^*\mathbb{R}$. For comparison, we can represent a complex number $a + bi$ by a two-tuple $\langle a, b\rangle$, a set-theoretical representation of which is
\begin{equation}
\langle a, b\rangle := \{a, \{a, b\}\}
\end{equation}
But using $\{a, \{a, b\}\} $ as set-theoretical representation for a single complex number does \textbf{absolutely not} mean that we have now identified the binary operation `addition' on the complex numbers with constructing a pair set, as in
\begin{equation}
z_1 + z_2 \equiv \{z_1, \{z_1, z_2\}\}
\end{equation}
So regarding addition of expanded reals, we merely have that e.g. the image under the operation addition of a Dedekind cut designated by a symbol `$a$' as in Eq. (\ref{eq:7}) and a countable set of two-tuples designated by a symbol `$b\omega$' as in Eq. (\ref{eq:8}) is a union of sets designated by a symbol `$b\omega + a$' as in Eq. (\ref{eq:8a}). An equation `$a$' + `$b\omega$' = `$b\omega + a$', with added quotation marks for clarity, is then a notation for $((a,b\omega),b\omega+a)\in +$. \hfill $\Box$
\end{Remark}

\noindent The following two definitions are also useful.
\begin{Definition}\
Let $x = \sum_{j=1}^n a_j\omega^{p_j}$ be an expanded real number. The \textbf{real part} of $x$ is then the number $Re(x)$ given by
\begin{equation}
\left\{
\begin{array}{l}
p_n = 0 \Rightarrow Re(x) = a_n \\
p_n > 0 \Rightarrow Re(x) = 0
\end{array}
\right.
\end{equation}
Likewise, the \textbf{hyperreal part} of $x$ is then the number $Hy(x)$ given by
\begin{equation}
\left\{
\begin{array}{l}
p_n = 0 \Rightarrow Hy(x) = x - a_n \\
p_n > 0 \Rightarrow Hy(x) = x
\end{array}
\right.
\end{equation}
So, for any expanded real number $x$ we have $x = Re(x) + Hy(x)$. \hfill $\Box$
\end{Definition}
\begin{Definition}\label{def:ReHy}\
Let $f: \mathbb{R} \rightarrow \ ^*_+\mathbb{R}$. Then the \textbf{real part of} $f$ is the function $f_{Re}: \mathbb{R} \rightarrow \ ^*_+\mathbb{R}$ given by
\begin{equation}
f_{Re}: x \mapsto Re(f(x))
\end{equation}
Likewise the \textbf{hyperreal part of} $f$ is the function $f_{Hy}: \mathbb{R} \rightarrow \ ^*_+\mathbb{R}$ given by
\begin{equation}
f_{Hy}: x \mapsto Hy(f(x))
\end{equation}
$\Box$
\end{Definition}

\noindent The space of all functions $f: \mathbb{R} \rightarrow \ ^*_+\mathbb{R}$ then forms a vector algebra over $\mathbb{R}$, when function addition, scalar multiplication, and function multiplication are defined naturally, so
$(f+g)(x) = f(x) + g(x)$, $(\alpha\cdot f)(x) = \alpha f(x)$, and $(f\cdot g)(x) = f(x)g(x)$. In particular, for any $f: \mathbb{R} \rightarrow \ ^*_+\mathbb{R}$ we have $f = f_{Re} + f_{Hy}$.\\
\ \\
In this function space $^*_+\mathbb{R}^\mathbb{R}$ we now define expanded real delta functions $\alpha\delta_\beta$ as follows:
\begin{Definition}\label{def:2}\
For any $\alpha,\beta\in \mathbb{R}$, the \textbf{expanded real delta function} $\alpha\delta_\beta$ is a function $\alpha\delta_\beta: \mathbb{R} \rightarrow \ ^*_+\mathbb{R}$ given by
\begin{enumerate}[(i)]
\item $x \neq \beta \Rightarrow \alpha\delta_\beta(x) = 0$
\item $x = \beta \Rightarrow \alpha\delta_\beta(x) = \alpha\omega$
\end{enumerate}
Note that for $f = \alpha\delta_\beta$ we have $f = f_{Hy}$. \hfill $\Box$
\end{Definition}
\noindent We can now define the integral over $\mathbb{R}$ of an arbitrary function $f: \mathbb{R} \rightarrow \ ^*_+\mathbb{R}$:
\begin{Definition}\label{def:7}\
Let $\mathcal{R}^1(\mathbb{R})$ be the set of Riemann integrable functions on $\mathbb{R}$, and let the set of all integrable expanded real functions on $\mathbb{R}$ be denoted by $^*_+\mathcal{R}^1(\mathbb{R})$. Let $f \in \ ^*_+\mathbb{R}^\mathbb{R}$; then $f \in \ ^*_+\mathcal{R}^1(\mathbb{R})$ \textbf{if and only if}
\begin{equation}
f_{Re} \in \mathcal{R}^1(\mathbb{R})
\end{equation}
\begin{equation}\label{eq:sum}
f_{Hy} = \sum_{n=1}^\infty \alpha_n\delta_{\beta_n}
\end{equation}
for some convergent series $\sum_{n=1}^\infty \alpha_n = s \in \mathbb{R}$; if $\int_{-\infty}^{+\infty} f_{Re}(x)dx = h$ then
\begin{equation}
\int_{-\infty}^{+\infty}f(x) dx = \int_{-\infty}^{+\infty} f_{Re}(x)dx + \int_{-\infty}^{+\infty} f_{Hy}(x)dx := h  + s
\end{equation}
$\Box$
\end{Definition}

\begin{Remark}\rm\
Def. \ref{def:7} has, thus, to be taken as an \textbf{axiom} of expanded real integral calculus: this does not ``follow'' from ordinary real integral calculus, nor can it be derived from the axioms of the hyperreal number field.\hfill$\Box$
\end{Remark}

\noindent Def. \ref{def:7} thus says that an expanded real function $f$ on $\mathbb{R}$ is integrable if and only if the real part of $f$ is Riemann integrable and the hyperreal part of $f$ is a countable sum of expanded real delta functions with the coefficients forming a convergent series. A corollary of Def. \ref{def:7} is that for the integral of the expanded real delta function $\alpha\delta_\beta$ over $\mathbb{R}$ we have
\begin{equation}\label{eq:9}
\int_{-\infty}^{+\infty}\alpha\delta_\beta(x)dx := \alpha
\end{equation}
\noindent As a set, an expanded real delta function $\alpha\delta_\beta$ can simply be identified with its graph:
\begin{equation}
\alpha\delta_\beta := \{ \langle x, \xi \rangle \in \mathbb{R}\times\ ^*_+\mathbb{R} \ | \ x \neq \beta \Rightarrow \xi = 0 \wedge x = \beta \Rightarrow \xi = \alpha\omega\}
\end{equation}
Each expanded real delta function $\alpha\delta_\beta$ is thus a special function with domain $\mathbb{R}$ that satisfies both clauses of Eq. (\ref{eq:0}). The function $\delta_0 = 1\delta_0$ in particular is thus a function on the reals with the properties of the Dirac delta $\delta$ displayed in Eq. (\ref{eq:1})---the dismissal of the Dirac delta by Von Neumann was thus premature!

\section{Generalization and conclusion}
Def. \ref{def:7} can be generalized to integrable expanded real functions on $\mathbb{R}^n$. For that matter, we need the following definition:
\begin{Definition}\label{def:permutation}\ Let $S_n$ be the group of permutations on $n$ letters. Let $\sigma \in S_n$. Then a \textbf{permutation of variables} $\pi_\sigma:\ ^*_+\mathbb{R}^{\mathbb{R}^n} \rightarrow\ ^*_+\mathbb{R}^{\mathbb{R}^n}$ is a function on the space of all functions from $\mathbb{R}^n$ to $\ ^*_+\mathbb{R}$, such that the image $\pi_\sigma(f)$ of an arbitrary function $f: \mathbb{R}^n \rightarrow\ ^*_+\mathbb{R}$ under $\pi_\sigma$ is given by
\begin{equation}\label{eq:10}
\pi_\sigma(f): (x^1,x^2, \ldots, x^n) \mapsto f(x^{\sigma(1)}, x^{\sigma(2)}, \ldots, x^{\sigma(n)})
\end{equation}
$\Box$
\end{Definition}

\begin{Example}\rm Let $\sigma \in S_3$ be the permutation $\left(\begin{array}{ccc}1 & 2 & 3 \\2 & 3 & 1 \\\end{array}\right)$, and let the function $f: \mathbb{R}^3 \rightarrow\ ^*_+\mathbb{R}$ be given by
$f: (x^1,x^2,x^3)\mapsto \sin (x^1x^2) \ \delta_0(x^3)$. Then $\pi_{\sigma}$ is a permutation of variables such that
\begin{equation}\label{eq:11}
\pi_{\sigma}(f): (x^1,x^2,x^3) \mapsto \sin (x^2x^3) \delta_0(x^1)
\end{equation}
$\Box$
\end{Example}
\noindent Furthermore, we extend Def. \ref{def:ReHy} to expanded real functions on $\mathbb{R}^n$. That is, for any function $f:\mathbb{R}^n\rightarrow\ ^*_+\mathbb{R}$ the real and hyperreal parts of $f$ are the functions $f_{Re}$ and $f_{Hy}$ given by
\begin{equation}
f_{Re}: (x^1, \ldots, x^n) \mapsto Re(f(x^1, \ldots, x^n))
\end{equation}
\begin{equation}
f_{Hy}: (x^1, \ldots, x^n) \mapsto Hy(f(x^1, \ldots, x^n))
\end{equation}
(Def. \ref{def:ReHy} can, of course, be extended to expanded real functions on any non-empty set $X$.)
\begin{Definition}\label{def:12}\
Let $\mathcal{R}^1(\mathbb{R}^n)$ be the set of Riemann integrable functions on $\mathbb{R}^n$, and let the set of all integrable expanded real functions on $\mathbb{R}^n$ be denoted by $^*_+\mathcal{R}^1(\mathbb{R}^n)$.
Let $f:\mathbb{R}^n\rightarrow\ ^*_+\mathbb{R}$; then $f \in \ ^*_+\mathcal{R}^1(\mathbb{R}^n)$ \textbf{if and only if}
\begin{equation}
f_{Re} \in \mathcal{R}^1(\mathbb{R}^n)
\end{equation}
\begin{equation}\label{eq:sum2}
f_{Hy} = \sum_{m=1}^\infty g_m
\end{equation}
for expanded real functions $g_m$, each of which is the image of a tensor product of an integrable function $u_m \in \ ^*_+\mathcal{R}^1(\mathbb{R}^{n-1})$ and an expanded real delta function $\alpha_m\delta_{\beta_m} \in\ ^*_+\mathcal{R}^1(\mathbb{R})$ under a permutation of variables $\pi_{\sigma_m}$ such that the integrals of $g_m$ and $f_{Hy}$ are finite:
\begin{equation}
g_m = \pi_{\sigma_m}(u_m \otimes \alpha_m\delta_{\beta_m})
\end{equation}
\begin{equation}
u_m \otimes \alpha_m\delta_{\beta_m}: (x^1, x^2, \ldots, x^n) \mapsto u_m(x^1, x^2, \ldots, x^{n-1})\cdot \alpha_m\delta_{\beta_m} (x^n)
\end{equation}
\begin{equation}
\pi_{\sigma_m}(u_m \otimes \alpha_m\delta_{\beta_m}):(x^1, x^2, \ldots, x^n) \mapsto u_m(x^{\sigma_m(1)}, \ldots, x^{\sigma_m(n-1)}) \alpha_m\delta_{\beta_m}(x^{\sigma_m(n)})
\end{equation}
\begin{equation}
\int_{-\infty}^\infty\cdots\int_{-\infty}^\infty u_m(x^1, \ldots, x^{n-1})dx^1\ldots dx^{n-1}\int_{-\infty}^\infty \alpha_m\delta_{\beta_m} (x^n)dx^n = s_m
\end{equation}
\begin{equation}
\int_{-\infty}^\infty\cdots\int_{-\infty}^\infty f_{Hy}(x^1, \ldots, x^n)dx^1\ldots dx^n = \sum_{m=1}^\infty s_m \in \mathbb{R}
\end{equation}
$\Box$
\end{Definition}

\begin{Example}\label{ex:1}\rm
Let $n=2$. Let the function $\alpha\delta^2_{(\beta_1, \beta_2)}:\mathbb{R}^2 \rightarrow \ ^*_+\mathbb{R}$ be given by
\begin{equation}
\alpha\delta^2_{(\beta_1, \beta_2)}: (x^1,x^2) \mapsto \alpha\delta_{\beta_1}(x^1)\delta_{\beta_2}(x^2)
\end{equation}
Then $\alpha\delta^2_{(\beta_1, \beta_2)}$ can be written as a tensor product:
\begin{equation}\label{eq:f}
\alpha\delta^2_{(\beta_1, \beta_2)} = \alpha\delta_{\beta_1}\otimes1\delta_{\beta_2}
\end{equation}
For the integral of $\alpha\delta^2_{(\beta_1, \beta_2)}$ over $\mathbb{R}^2$ we then have
\begin{equation}\label{eq:last}
\int_{-\infty}^{\infty}\int_{-\infty}^{\infty}\alpha\delta^2_{(\beta_1, \beta_2)}(x^1, x^2)dx^1dx^2 =
\int_{-\infty}^{\infty}\alpha\delta^2_{\beta_1}(x^1)dx^1\int_{-\infty}^{\infty}1\delta_{\beta_2}(x^2)dx^2 = \alpha
\end{equation}
$\Box$
\end{Example}
\noindent Expanded real functions $\alpha\delta^n_{(\beta_1, \ldots, \beta_n)}:\mathbb{R}^n \rightarrow \ ^*_+\mathbb{R}$, which can be written as a tensor product of $n$ expanded real delta functions just like the function $\alpha\delta^2_{(\beta_1, \beta_2)}$ in Eq. (\ref{eq:f}), can henceforth be called \textbf{Dirac delta functions}---as opposed to the Dirac delta \textbf{distributions}, defined according to Eq. (\ref{eq:abusive}).\\
\ \\
Concluding, this paper has introduced \emph{Dirac delta functions} $\alpha\delta_\beta$ as special hyperreal-valued functions of a real variable. If we interpret the function value as a density, we can speak of \emph{infinitely high densities} in those points where the value is nonzero, since $|\alpha\cdot\omega|=\infty$ for any $\alpha \in \mathbb{R}$. One might object that these new functions do not enable one to do previously impossible calculations or to prove new theorems in pure mathematics. That is true, but for the mathematical modeling of physical systems the newly introduced functions are nevertheless advantageous compared to the existing objects that capture the idea of the Dirac delta: unlike these existing objects, the new functions enable one to model the state of a physical system in which infinitely high densities occur---such as the mass density in a system made up of point particles---by representing the distribution of the density over space with a function on $\mathbb{R}^n$. In particular, these functions will be applied in the construction of a Planck scale model of the Elementary Process Theory, an  abstract scheme with a physical interpretation introduced in \cite{Cabbolet}.



\begin{thebibliography}{7}

\bibitem{Dirac}
P.A.M. Dirac, \emph{Proc. Roy. Soc. A} \textbf{113}(765), 621--641 (1927)

\bibitem{VonNeumann}
J. Von Neumann, \emph{Mathematische Grundlagen der Quantenmechanik}, Springer Verlag, Berlin, p. 14 (1932)

\bibitem{Schwartz}
L. Schwartz, \emph{Theorie des Distributions. Tom I}, Hermann, Berlin (1950)

\bibitem{Rudin}
W. Rudin, \emph{Real and Complex Analysis}, McGraw-Hill, New York (1966)

\bibitem{Robinson}
A. Robinson, \emph{Non-standard Analysis}, Princeton University Press, Princeton (1966)

\bibitem{Todorov}
T. Todorov, \emph{Proc. Am. Math. Soc.} \textbf{110}(4), 1143--1144 (1990)

\bibitem{Dannon}
H.V. Dannon, \emph{Gauge Institute Journal of Math and Physics} \textbf{8}(1), 1--49 (2012)

\bibitem{Cabbolet}
M.J.T.F. Cabbolet, \emph{Ann. Phys.} \textbf{522}, 699-738 (2010); \textbf{523}, 990-994 (2011); \textbf{528}, 626-627 (2016)
\end{thebibliography}
\end{document}